\documentclass[oneside]{amsart}

\usepackage[letterpaper,body={14.6cm,23.5cm}, mag=1000]{geometry}
\usepackage{amssymb}
\usepackage{amsthm}
\usepackage{amscd}

\numberwithin{equation}{section}
\theoremstyle{plain}
\newtheorem{thm}[equation]{Theorem}
\newtheorem{cor}[equation]{Corollary}
\newtheorem{lemma}[equation]{Lemma}

\newtheorem{prop}[equation]{Proposition}

\newtheorem*{thma}{Theorem A}
\newtheorem*{thmb}{Theorem B}
\newtheorem*{thmb'}{Theorem B$^\prime$}
\theoremstyle{definition}

\newcommand{\dlabel}[1]{\ifmmode \text{\ttfamily \upshape [#1] } \else
{\ttfamily \upshape [#1] }\fi \label{#1}}

\newcommand{\C}{\operatorname{C} }

\newcommand{\Z}{\operatorname{Z} }

\newcommand{\gen}[1]{\left < #1 \right >}
\newcommand{\Aut}{\operatorname{Aut} }

\newcommand{\Inn}{\operatorname{Inn} }

\begin{document}

\baselineskip 15pt

\title{Central quotient versus commutator subgroup of groups}

\author{Manoj K.~Yadav}

\address{School of Mathematics, Harish-Chandra Research Institute \\
Chhatnag Road, Jhunsi, Allahabad - 211 019, INDIA}

\email{myadav@hri.res.in}


\begin{abstract}
In 1904, Issai Schur proved the following result. If $G$ is an arbitrary group such that $G/\Z(G)$ is finite, where $\Z(G)$ denotes the center of the group $G$, then the commutator subgroup of $G$ is finite. A partial converse of this result was proved by B. H. Neumann in 1951. He proved that if  $G$ is a finitely generated group with finite commutator subgroup, then $G/\Z(G)$ is finite. In this short note,  we exhibit  few arguments of  Neumann, which provide further generalizations of converse of  the above mentioned result of Schur. We classify all finite groups $G$ such that $|G/\Z(G)| = |\gamma_2(G)|^d$, where $d$ denotes the number of elements in a minimal generating set for $G/\Z(G)$.   Some problems and  questions are posed in the sequel. 
\end{abstract}
\maketitle

\section{Introduction}

In 1951,  Neumann \cite[Theorem 5.3]{bN51} proved the following result: \emph{If the index of $\Z(G)$ in $G$ is finite, then $\gamma_2(G)$ is finite, where $\Z(G)$ and $\gamma_2(G)$ denote the center and the commutator subgroup of $G$ respectively.} He mentioned \cite[End of page 237]{bN54} that this result can be obtained from an implicit idea of  Schur \cite{iS04}, and his proof also used Schur's basic idea. However there is no mention of this fact in \cite{bN51} in which Schur's paper is also cited.  In this note, this result will be termed as `the Schur's theorem'. Neumann  also provided a partial converse of  the Schur's theorem \cite[Corollary 5.41]{bN51} as follows: {\it  If   $G$ is finitely generated by $k$ elements and $\gamma_2(G)$ is finite, then $G/\Z(G)$ is  finite, and bounded by  $|G/\Z(G)| \le |\gamma_2(G)|^k$}. 

Our first motivation of writing this note is to exhibit an idea  of Neumann\cite[page 179]{bN51} which proves much more than what is said above on  converse of the Schur's theorem.  We quote the text here (with a minor modification in the notations): 

\emph{``Let $G$ be generated by $g_1, g_2, \ldots, g_k$. Then	
\[\Z(G) = \cap_{\kappa=1}^{\kappa=k} C_G(g_{\kappa});\]
for, an element of $G$ lies in the centre if and only if it (is permutable) commutes with all the generators of $G$. If $G$ is an FC-group (group whose all conjugacy classes are of finite length), then $|G:C_G(g_{\kappa})|$ is finite for $1 \le \kappa \le k$, and $\Z(G)$, as intersection of a finite set of subgroups of finite index, also has finite index. The index of the intersection of two subgroups does not exceed the product of the indices of the subgroups: hence in this case one obtains an upper bound for the index of the centre, namely
\[|G : \Z(G)| \le  \Pi_{\kappa=1}^{\kappa=k} |G : C_G(g_{\kappa})|."\]}

Just a soft staring at the quoted text for a moment or two  suggests the following. The conclusion does not require the group $G$ to be FC-group. It only  requires the finiteness of the conjugacy classes of the generating elements.  If a generator of the group $G$ is contained in $\Z(G)$, one really does not need to count it. Thus the argument works perfactly well even if $G$ is generated by infinite number of elements, all but finite  of them  lie in the  center of $G$.   Thus the following result holds true.

\begin{thma}
 Let $G$ be an arbitrary group such that $G/\Z(G)$ is finitely generated by $x_1\Z(G),$ $x_2\Z(G),$ $\ldots, x_t\Z(G)$ and the conjugacy class of $x_i$ in $G$ is of finite length for $1 \le i \le t$. Then $G/\Z(G)$ is finite. Moreover $|G/\Z(G)| \le \Pi_{i = 1}^t |x_i^G|$ and $\gamma_2(G)$ is finite, where $x_i^G$ denotes the conjugacy class of $x_i$ in $G$.
\end{thma}

Neumann's result \cite[Corollary 5.41]{bN51} was reproduced by  Hilton \cite[Theorem 1]{pH01}. It seems  that Hilton was not aware of  Neumann's result.  This lead two more publications \cite{pN10} and \cite{bS10} dedicated to proving special cases of Theorem A.

Converse of the Schur's theorem is not true in general as shown by infinite extraspecial $p$-groups, where $p$ is an  odd prime. It is interesting to know that example of such a $2$-group also exists, which is mentioned on page 238 (second para of Section 3) of \cite{bN54}. It is a central product of infinite copies of quaternion groups of order $8$ amalgamated at the center of order $2$.

Our second motivation of  writing this note is to provide a modification of  an innocent looking result of Neumann \cite[Lemma 2]{bN76}, which allows us to say little more on converse of the Schur's theorem.   A modified version of this lemma is the following. 

\begin{lemma}\label{BHNlemma}
Let $G$ be an arbitrary group having a normal abelian subgroup $A$ such that the index of $C_G(A)$ in $G$ is finite and  $G/A$ is finitely generated by $g_1A, g_2A, \ldots, g_rA$, where  $|g_i^G| < \infty$ for  $1 \le i \le r$. Then $G/Z(G)$ is finite.
\end{lemma} 

This lemma helps proving the first three statements of the following result.

\begin{thmb}
For an arbitrary group $G$, $G/\Z(G)$ is finite if any one of the following holds true:

(i) $\Z_2(G)/\Z(\Z_2(G))$ is finitely generated and $\gamma_2(G)$ is finite.

(ii) $G/\Z(\Z_2(G))$ is finitely generated and $G/(\Z_2(G)\gamma_2(G))$ is finite.

(iii) $\gamma_2(G)$ is finite and $\Z_2(G) \le \gamma_2(G)$.

(iv) $\gamma_2(G)$ is finite  and $G/\Z(G)$ is purely non-abelian.
\end{thmb}

Our final motivation is to provide a classification of all  groups $G$ upto isoclinism (see Section 3 for the definition) such that $|G/\Z(G)| = |\gamma_2(G)|^d$ is finite, where $d$ denotes the number of elements in a minimal generating set for $G/\Z(G)$, discuss example in various situations and pose some problems.
 We conclude this section with fixing some notations. For  an arbitrary group $G$,  by $\Z(G)$, $\Z_2(G)$ and $\gamma_2(G)$ we denote the center, the second center and the commutator subgroup of $G$ respectively.  For $x \in G$, $[x,G]$ denotes the set $\{[x,g] \mid g \in G\}$. Notice that $|[x,G]| = |x^G|$, where $x^G$ denotes the conjugacy class of $x$ in $G$. If $[x,G] \subseteq \Z(G)$, then $[x,G]$ becomes a subgroup of $G$.  For a subgroup $H$ of $G$, $\C_{G}(H)$ denotes the centralizer of $H$ in $G$ and for an element $x \in G$, $\C_{G}(x)$ denotes the centralizer of $x$ in $G$.

\section{Proofs}

We start with the proof of Lemma \ref{BHNlemma}, which  is essentialy same as the one given by  Neumann.  

\noindent {\bf Proof of Lemma \ref{BHNlemma}.}  Let $G/A$ be generated by $g_1A, g_2A, \ldots, g_rA$ for some $g_i \in G$, where $1 \le i \le r < \infty$. Let $X := \{g_1, g_2, \ldots, g_r\}$ and $A$ be generated by a set $Y$. Then $G = \gen{X \cup Y}$ and $\Z(G) = \C_G(X) \cap C_G(Y)$. Notice that $C_G(A) = C_G(Y)$. Since $C_G(A)$ is of finite index,  $C_G(Y)$ is also of finite index in $G$. Also, since $|g_i^G|\ < \infty$ for $1 \le i \le r$,  $C_G(X)$  is of finite index in $G$. Hence the index of $Z(G)$ in $G$ is finite and the proof is complete. \hfill $\Box$

\vspace{.1in}

Proof of Theorem A can be made quite precise by using Lemma  \ref{BHNlemma}.

\vspace{.1in}

\noindent {\bf Proof of Theorem A.}  Taking $A = \Z(G)$ in Lemma \ref{BHNlemma}, it follows that $G/\Z(G)$ is finite. Moreover,  \[|G/\Z(G)| = |G/\cap_{i=1}^tC_G(x_i)| \le \Pi_{i=1}^t |G:C_G(x_i)| = \Pi_{i=1}^t |[x_i, G]| = \Pi_{i = 1}^t |x_i^G|.\]   
That $\gamma_2(G)$ is finite now follows from the Schur's theorem. \hfill $\Box$

\vspace{.1in}

For the proof of Theorem B we need the following result of  Hall \cite{pH56} and the subsequent proposition.

\begin{thm}\label{thm1}
If $G$ is an arbitrary group such that $\gamma_2(G)$ is finite, then $G/\Z_2(G)$ is finite.
\end{thm}

Explicit bounds on the order of $G/\Z_2(G)$ were first given by  Macdonald \cite[Theorem 6.2]{iM61} and later on improved by   Podoski and  Szegedy \cite{PS05} by showing that if $|\gamma_2(G)/(\gamma_2(G) \cap \Z(G))| = n$, then $|G/\Z_2(G)| \le n^{c \;log_2 n}$ with $c = 2$.

\begin{prop}\label{prop}
Let $G$ be an arbitrary group such that $\gamma_2(G)$ is finite and  $G/\Z(G)$ is infinite. Then $G/\Z(G)$ has an infinite abelian group as a direct factor.
\end{prop}
\begin{proof}
Since $\gamma_2(G)$ is finite, by Theorem \ref{thm1} $G/\Z_2(G)$ is finite. Thus $\Z_2(G)/\Z(G)$ is infinite. Again using the finiteness of $\gamma_2(G)$, it follows that the exponent of $\Z_2(G)/\Z(G)$ is finite. Hence by  \cite[Theorem 17.2]{lF70} $\Z_2(G)/\Z(G)$ is a direct sum of cyclic groups. Let $G/\Z_2(G)$ be generated by $x_1\Z_2(G), \ldots, x_r\Z_2(G)$ and $H := \gen{x_1, \ldots, x_r}$. Then it follows that modulo $\Z(G)$, $H \cap \Z_2(G)$ is finite. Thus we can write 
\[\Z_2(G)/\Z(G) = \gen{y_1\Z(G)} \times \cdots \times \gen{y_s\Z(G)} \times \gen{y_{s+1}\Z(G)} \times \cdots,\]
such that $(H \cap \Z_2(G))\Z(G)/\Z(G) \le \gen{y_1\Z(G)} \times \cdots \times \gen{y_s\Z(G)}$. It now follows that the infinite abelian group $\gen{y_{s+1}\Z(G)} \times \cdots$ is a direct factor of $G/\Z(G)$, and  the proof is complete. \hfill $\Box$

\end{proof}

We are now ready to prove Theorem B.
\vspace{.1in}

\noindent {\bf Proof of Theorem B.}  Since $\gamma_2(G)$ is finite, it follows from Theorem \ref{thm1}  that $G/\Z_2(G)$ is finite. Now using the fact that $\Z_2(G)/\Z(\Z_2(G))$ is finitely generated, it follows that $G/\Z(\Z_2(G))$ is finitely generated. Take $\Z(\Z_2(G)) = A$. Then notice that $A$ is a normal abelian subgroup of $G$ such that the index of $C_G(A)$ in $G$ is finite, since $\Z_2(G) \le  C_G(A)$. Hence by Lemma \ref{BHNlemma}, $G/\Z(G)$ is finite, which proves (i).

Again take $\Z(\Z_2(G)) = A$ and notice that $\Z_2(G)\gamma_2(G) \le \C_G(A)$. (ii) now directly follows from Lemma \ref{BHNlemma}.  If $\Z_2(G) \le \gamma_2(G)$, then $\Z_2(G)$ is abelian. Thus (iii) follows from (i).  Finally, (iv) follows from Proposition \ref{prop}. This completes the proof of the theorem. \hfill $\Box$

We conclude this section with an extension of Theorem A in terms of conjugacy class-preserving automorphisms of given group $G$. An automorphism $\alpha$ of an arbitrary group $G$ is called \emph{(conjugacy) class-preserving} if $\alpha(g) \in g^G$ for all $g \in G$. We denote the group of all class-preserving automorphisms of $G$ by $\Aut_c(G)$. Notice that $\Inn(G)$, the group of all inner automorphisms of $G$, is a normal subgroup of $\Aut_c(G)$ and $\Aut_c(G)$ acts trivially on the center of $G$.
A detailed survey on class-preserving automorphisms of finite $p$-groups can be found in \cite{mY11}.

Let $G$ be the group as in the statement of Theorem A. Then $G$ is generated by $x_1, x_2, \ldots, x_t$ along with $\Z(G)$.  Since $\Aut_c(G)$ acts trivially on the center of $G$, it follows that 
\begin{equation}\label{eqn1}
|\Aut_c(G)| \le \Pi_{i=1}^t |x_i^G|
\end{equation}
 as there are only $|x_i^G|$ choices for the image of each $x_i$ under any class-preserving automorphism. 
Since $|x_i^G|$ is finite for each $x_i$, $1 \le i \le t$, it follows that  $|\Aut_c(G)| \le \Pi_{i=1}^t |x_i^G|$ is finite.

We have proved the following result of which Theorem A is a corollary, because $|G/\Z(G)| = |\Inn(G)| \le |\Aut_c(G)|$.

\begin{thm}
 Let $G$ be an arbitrary group such that $G/\Z(G)$ is finitely generated by $x_1\Z(G),$ $x_2\Z(G), $ $\ldots, x_t\Z(G)$ and the conjugacy class of $x_i$ in $G$ is of finite length for $1 \le i \le t$. Then $\Aut_c(G)$ is finite. Moreover $|\Aut_c(G)| \le \Pi_{i = 1}^t |x_i^G|$ and $\gamma_2(G)$ is finite.
\end{thm}

Proof of Theorem A is also reproduced using $IA$-automorphisms (automorphisms of a group that induce identity on the abelianization) in \cite[Theorem 2.3]{GK13}. Proof goes on the same way as in the case of  class-preserving automorphisms.

\section{Groups with maximal central quotient}

We start with the following concept due to  Hall \cite{pH40}. For a group $X$, the commutator map $a_X: X/\Z(X) \times X/\Z(X) \to \gamma_2(X)$ given by $a_X(x_1\Z(X), x_2\Z(X)) = [x_1, x_2]$ is well defined.
Two groups $K$ and $H$ are said to be \emph{isoclinic} if 
there exists an  isomorphism $\phi$ of the factor group
$\bar K = K/\Z(K)$ onto $\bar{H} = H/\Z(H)$, and an isomorphism $\theta$ of
the subgroup $\gamma_{2}(K)$ onto  $\gamma_{2}(H)$
such that the following diagram is commutative
\[
 \begin{CD}
   \bar K \times \bar K  @>a_G>> \gamma_{2}(K)\\
   @V{\phi\times\phi}VV        @VV{\theta}V\\
   \bar H \times \bar H @>a_H>> \gamma_{2}(H).
  \end{CD}
\]
The resulting pair $(\phi, \theta)$ is called an \emph{isoclinism} of $K$ 
onto $H$. Notice that isoclinism is an equivalence relation among groups.

 The following proposition (also see  Macdonald's result \cite[Lemma 2.1]{iM61}) is important for the rest of this section.

\begin{prop}\label{3prop1}
 Let $G$ be a group such that $G/\Z(G)$ is finite. Then there exists a finite group $H$  isoclinic to the group $G$ such that $\Z(H) \le \gamma_2(H)$. Moreover if $G$ is a $p$-group, then $H$ is also a $p$-group.
\end{prop}
\begin{proof}
 Let $G$ be the given group. Then by Schur's theorem $\gamma_2(G)$ is finite. Now it follows from a result of  Hall \cite{pH40} that there exists a group $H$ which is isoclinic to $G$ and $\Z(H) \le \gamma_2(H)$. Since $|\gamma_2(G)| = |\gamma_2(H)|$ is finite, $\Z(H)$ is finite. Hence, by the definition of isoclinism, $H$ is finite. Now suppose that $G$ is a $p$-groups. Then it follows that $H/\Z(H)$ as well as $\gamma_2(H)$ are $p$-groups. Since $\Z(H) \le \gamma_2(H)$, this implies that $H$ is a $p$-group.    \hfill $\Box$

\end{proof}

For an arbitrary group $G$ with finite $G/\Z(G)$, we have 
\begin{equation}\label{3eqn1}
|G/\Z(G)| \le |\gamma_2(G)|^d,
\end{equation}
where $d = d(G/\Z(G))$.
For simplicity we say that a group $G$ has \emph{Property A} if $G/\Z(G)$ is finite and equality holds in \eqref{3eqn1} for $G$. We are now going to classify, upto isoclinism,  all groups $G$ having Property A.

Let $G$ be an arbitrary group having Property A. Then by Proposition \ref{3prop1} there exists a finite group $H$ isoclinic to $G$ and, by the definition of isoclinism, H has Property A. Thus for classifying all groups $G$, upto isoclinism,  having Property A, it is sufficient to classify all finite group with this property.  

Let us first consider  non-nilpotent finite groups. For such group we prove the following result in \cite{DPY} 

\begin{thm}\label{thmDPY}
There is no non-nilpotent  group $G$ having Property A.
\end{thm}

So we only need to consider finite nilpotent groups. Since a finite nilpotent group is a direct product of it's Sylow $p$-subgroups, it is sufficient to classify finite $p$-groups admitting Property A. 
Obviously, all abelian groups admit Property A.   Perhaps the simplest examples of non-abelian groups having Property A are finite extraspecial $p$-groups.   The class of $2$-generated finite capable nilpotent groups with cyclic commutator subgroup also admits Property A. A group $G$ is said to be \emph{capable} if there exists a group $H$ such that $G \cong H/\Z(H)$.  Isaacs \cite[Theorem 2]{mI01} proved: \emph{Let $G$ be finite and capable, and suppose that $\gamma_2(G)$ is cyclic and that all elements of order $4$ in $\gamma_2(G)$ are central in $G$. Then $|G/\Z(G)| \le |\gamma_2(G)|^2$, and equality holds if $G$ is nilpotent}. So $G$ admits Property A, if $G$ is a nilpotent group as in this statement and $G$ is $2$-generated.  A complete classification of $2$-generated finite capable $p$-groups of class $2$ is given in  \cite{MM10}. 

Motivated by finite extraspecial $p$-groups, a more general class of groups $G$ admitting Property A can be constructed as follows. For any positive integer $m$, let $G_1, G_2, \ldots, G_m$ be  $2$-generated finite $p$-groups such that $\gamma_2(G_i) = \Z(G_i) \cong X$ (say) is cyclic of order  $q$ for $1 \le i \le m$, where $q$ is some power of $p$. Consider the central product 
\begin{equation}\label{ygroup}
Y = G_1 *_X G_2 *_X \cdots *_X G_m
\end{equation}
 of $G_1, G_2, \ldots, G_m$ amalgamated at $X$ (isomorphic to cyclic commutator subgroups $\gamma_2(G_i)$, $1 \le i \le m$).  Then $|Y| = q^{2m+1}$ and $|Y/\Z(Y)| = q^{2m} = |\gamma_2(Y)|^{d(Y)}$, where  $d(Y) = 2m$ is the  number of elements in any minimal generating set for $Y$. Thus  $Y$ has Property A. Notice that in all of the above examples, the commutator subgroup is cyclic. Infinite groups having Property A can be easily obtained by taking a direct product of an infinite abelian group with any finite group having Property A.

We now proceed to showing  that any finite $p$-group $G$ of class $2$ having Property A is isoclinic to a group $Y$ defined in \eqref{ygroup}.  

Let $x \in \Z_2(G)$ for a group $G$. Then, notice that $[x, G]$ is a central subgroup of $G$. We have the following simple but useful result.

\begin{lemma}\label{3lemma1}
  Let $G$ be an arbitrary group such that $\Z_2(G)/\Z(G)$ is finitely generated by $x_1\Z(G),$ $ x_2\Z(G), \ldots, x_t\Z(G)$ such that $exp([x_i, G])$ is finite for $1 \le i \le t$. Then 
\[|\Z_2(G)/\Z(G)| = \prod_{i=1}^t exp([x_i,G]).\]
\end{lemma}
\begin{proof}
By the given hypothesis $exp([x_i,G])$ is finite for all $i$ such that $1 \le i \le t$. Suppose that $exp([x_i,G]) = n_i$. Since $[x_i, G] \subseteq \Z(G)$, it follows that  $[x_i^{n_i}, G] = [x_i, G]^{n_i} = 1$. Thus $x_i^{n_i} \in \Z(G)$ and no smaller power of $x_i$ than $n_i$ can lie in $\Z(G)$, which implies that the order of $x_i\Z(G)$ is $n_i$.  Since $\Z_2(G)/\Z(G)$ is abelian, we have $|\Z_2(G)/\Z(G)| = \prod_{i=1}^t exp([x_i,G])$.  \hfill $\Box$

\end{proof}

Let $\Phi(X)$ denote the Frattini subgroup of  a group $X$. The following result provides some structural information of $p$-groups of class $2$  admitting Property A.

\begin{prop}\label{3prop2}
Let $H$ be a finite  $p$-group of class $2$ having Property A and  $\Z(H) = \gamma_2(H)$. Then 

(i) $\gamma_2(H)$ is cyclic;

(ii) $H/\Z(H)$ is homocyclic;

(iii)  $[x, H] = \gamma_2(H)$ for all $x \in H - \Phi(H)$;

(iv) $H$ is minimally generated by even number of elements.
\end{prop}
\begin{proof}
Let $H$ be the group given in the statement,  which is minimally generated by $d$ elements $x_1, x_2 \ldots, x_d$ (say). Since $\Z(H) = \gamma_2(H)$, it follows that $H/\Z(H)$ is minimally generated by  $x_1\Z(H),$ $ x_2\Z(H), \ldots, x_d\Z(H)$.    Thus by the identity $|H/\Z(H)| = |\gamma_2(H)|^d$, it follows that order of $x_i\Z(H)$ is equal to $|\gamma_2(H)|$ for all $1 \le i \le d$. Since the exponent of $H/\Z(H)$ is equal to the exponent of $\gamma_2(H)$, we have that $\gamma_2(H)$ is cyclic and $H/\Z(H)$ is homocyclic. Now by Lemma \ref{3lemma1}, $|\gamma_2(H)|^d = |H/\Z(H)| = \prod_{i=1}^t exp([x_i,H])$. Since $[x_i, H] \subseteq \gamma_2(H)$,  this implies that $[x_i,H] = \gamma_2(H)$ for each $i$ such that $1 \le i \le d$.   Let $x$   be an arbitrary element in $H - \Phi(H)$. Then the set $\{x\}$ can always be extended to a minimal generating set of $H$. Thus it follows that $[x, H] = \gamma_2(H)$ for all $x \in H - \Phi(H)$.  This proves first three assertions. 

For the proof of (iv), we consider the group ${\bar H} = H/\Phi(\gamma_2(H))$.  Notice that both $H$ as well as  ${\bar H}$ are minimallay generated by $d$ elements. Since $[x, H] =  \gamma_2(H)$ for all $x \in H - \Phi(H)$, it follows that for no $x \in H - \Phi(H)$, ${\bar x} \in \Z({\bar H})$, where ${\bar x} = x\Phi(\gamma_2(H)) \in {\bar H}$. Thus it follows that 
$\Z({\bar H}) \le \Phi({\bar H})$. Also, since $\gamma_2(H)$ is cyclic, $\gamma_2({\bar H})$ is cyclic of order $p$. Thus it follows that  ${\bar H}$ is isoclinic to a finite extraspecial $p$-group, and therefore it is minimally generated by even number of elements.  Hence $H$ is also  minimally generated by even number of elements.  This completes the proof of the proposition. \hfill $\Box$

\end{proof}

Using the definition of isoclinism, we have
\begin{cor}
Let $G$ be a finite  $p$-group  of class $2$ admitting Property A. Then $\gamma_2(G)$ is cyclic and  $G/\Z(G)$ is homocyclic.
\end{cor}

We need the following important result.

\begin{thm}[\cite{BBC}, Theorem 2.1]\label{3thm1}
Let $G$ be a finite $p$-group of nilpotency class $2$ with cyclic center.  Then $G$ is a central product either of  two generator subgroups with cyclic center  or   two generator subgroups with cyclic center and a cyclic subgroup.
\end{thm}

\begin{thm}\label{thmcl2b}
Let $G$ be a  finite $p$-group of class $2$ having Property A.  Then $G$  is isoclinic to the group $Y$, defined in \eqref{ygroup}, for  suitable positive integers $m$ and $n$.
\end{thm}
\begin{proof}
Let $G$ be a group as in the statement. Then by Proposition \ref{3prop1} there exists a finite $p$-group $H$ isoclinic to $G$ such that $\Z(H) = \gamma_2(H)$. Obviously $H$ also satisfies $|H/\Z(H)| = |\gamma_2(H)|^d$, where $d$ denotes the number of elements in any minimal generating set of $G/\Z(G)$.  Then by Proposition \ref{3prop2}, $\gamma_2(H) = \Z(H)$ is cyclic of order $q = p^n$ (say)  for some positive integer $n$, and $H/\Z(H)$ is homocyclic of exponent $q$ and is of order $q^{2m}$ for some positive integer $m$.   Since $\Z(H) = \gamma_2(H)$ is cyclic, it follows from Theorem \ref{3thm1} that $H$ is a central product of $2$-generated groups $H_1, H_2, \ldots, H_m$. It is easy to see that $\gamma_2(H_i) = \Z(H_i)$ for $1 \le i \le m$ and $|\gamma_2(H)| = q$. This completes the proof of the theorem. \hfill $\Box$

\end{proof}

We would like to remark that Theorem \ref{thmcl2b} is also obtained in \cite[Theorem 11.2]{mY13} as a consequence on study of class-preserving automorphisms of finite $p$-group. But we have presented here a direct proof. 

Now we classify finite $p$-groups of nilpoency class larger than $2$. Consider the metacylic groups
\begin{equation}\label{Intgrp1}
K := \gen{x, y \mid x^{p^{r+t}} = 1, y^{p^r} = x^{p^{r+s}}, [x, y] = x^{p^t}},
\end{equation}
where $1 \le t < r$ and $0 \le s \le t$ ($t \ge 2$ if $p = 2$) are non-negative integers. Notice that the nilpotency class of $K$ is at least $3$. Since $K$ is generated by $2$ elements, it follows from \eqref{eqn1} that $|\Aut_c(K)| \le |\gamma_2(K)|^2 = p^{2r}$. It is not so difficult  to  see that $|\Inn(K)| = |K/\Z(K)| = p^{2r}$. Since $\Inn(K) \le \Aut_c(K)$, it follows that $|\Aut_c(K)| = |\Inn(K)| = |\gamma_2(K)|^2 = |\gamma_2(K)|^{d(K)}$ (That $\Aut_c(G) = \Inn(G)$, is, in fact, true for all finite metacylic $p$-groups). Thus $K$ admits Property A.  Furthermore, if $H$ is any $2$-generator group isoclinic to $K$, then it follows that $H$ admits Property A. For a finite $p$-groups having Property A, there always exists a $p$-group $H$ isoclinic $G$ such that $|H/\Z(H)| = |\gamma_2(H)|^d$, where $d = d(H)$. The following theorem now classifies, upto isoclinism, all finite $p$-groups $G$ of nilpotency class larger than $2$ having Property A.

\begin{thm}[Theorem 11.3, \cite{mY13}]\label{s11thm2}
Let $G$ be a finite $p$-group  of nilpotency class at least $3$.  Then  the following hold true.

(i) If $|G/\Z(G)| = |\gamma_2(G)|^d$, where $d = d(G)$, then $d(G) = 2$;

(ii) If  $|\gamma_2(G)/\gamma_3(G)| > 2$, then $|G/\Z(G)| = |\gamma_2(G)|^d$  if and only if $G$ is a $2$-generator group with cyclic commutator subgroup.  Furthermore, $G$ is isoclinic to the group $K$ defined in \eqref{Intgrp1} for suitable parameters;

(iii) If $|\gamma_2(G)/\gamma_3(G)| = 2$, then  $|G/\Z(G)| = |\gamma_2(G)|^d$  if and only if $G$ is a $2$-generator $2$-group of nilpotency class $3$ with elementary abelian $\gamma_2(G)$ of order $4$.
\end{thm}

It is clear that the groups $G$ occuring in Theorem \ref{s11thm2}(iii) are isoclinic to certain groups of order $32$. Using Magma (or GAP), one can easily show that such groups of order 32 are SmallGroup(32,k) for $k =6,7,8$  in the small group library.

We conclude this section with providing some different type of bounds on the central quotient of a given group.
If $|\gamma_2(G)\Z(G)/\Z(G)| = n$ is finite for a group $G$, then it follows from \cite[Theorem 1]{PS05} that  $|G/\Z_2(G)| \le  n^{2\;log_2 n}$. Using this and Lemma \ref{3lemma1}  we can also 
 provide an upper bound on the size of $G/\Z(G)$ in terms of  $n$, the rank of $\Z_2(G)/\Z(G)$ and exponents of certain sets of commutators (here these sets are really subgroups of $G$) of coset  representatives of generators of $\Z_2(G)/\Z(G)$ with the elements of $G$. This is given in the following theorem.

\begin{thm}\label{3thm2}
 Let $G$ be an arbitrary group. Let  $|\gamma_2(G)\Z(G)/\Z(G)| = n$ is finite and $\Z_2(G)/\Z(G)$ is finitely generated by $x_1\Z(G), x_2\Z(G), \ldots, x_t\Z(G)$ such that $exp([x_i, G])$ is finite for $1 \le i \le t$. Then
\[|G/\Z(G)| \le n^{2\;log_2 n} {\prod}_{i=1}^t exp([x_i, G]).\]
\end{thm}

\section{Problems and Examples}

Theorem B provides some conditions on a group $G$ under which $G/\Z(G)$ becomes finite. It is interesting to solve
\vspace{.1in}

\noindent {\bf Problem 1.} Let $G$ be an arbitrary group. Provide a set $\mathcal{C}$ of  optimal conditions on  $G$ such that $G/\Z(G)$ is finite if and only if all conditions in $\mathcal{C}$ hold true. 

As we know that there is no finite non-nilpotent group $G$ admitting Property A. Since $\Inn(G) \le \Aut_c(G)$, it is interesting to consider
\vspace{.1in}

\noindent{\bf Problem 2.} Classify all non-nilpotent finite group $G$ such that $|\Aut_c(G)| = |\gamma_2(G)|^d$, where $d = d(G)$.
\vspace{.1in}

A much stronger result than Theorem \ref{thmDPY} is known in the case when the  Frattini subgroup of $G$ is trivial. This is given in the following theorem of Herzog, Kaplan and Lev \cite[Theorem A]{HKL08} (the same result is also proved independently by Halasi and Podoski in \cite[Theorem 1.1]{HP08}).

\begin{thm}
Let $G$ be any non-abelian group with trivial Frattini subgroup. Then $|G/\Z(G)| < |\gamma_2(G)|^2$.
\end{thm}

The following result with the assertion similar to the preceding theorem is due to Isaacs \cite{mI01}.

\begin{thm}
 If $G$ is a capable finite group with cyclic $\gamma_2(G)$ and 
 all elements of order $4$ in $\gamma_2(G)$ are central in $G$, then $|G : Z(G)| \le |\gamma_2(G)|^2$. Moreover, equality holds if $G$ is nilpotent.
\end{thm}

So, there do exist nilpotent groups with comparatively small central quotient. A natural problem is the following.
\vspace{.1in}

\noindent{\bf Problem 3.} Classify all finite $p$-groups $G$ such that $|G : Z(G)| \le  |\gamma_2(G)|^2$.
\vspace{.1in}

Let $G$ be a finite nilpotent group of class $2$ minimally generated by $d$ elements. Then it follows from Lemma \ref{3lemma1} that $|G/\Z(G)| \le {\prod}_{i=1}^d exp([x_i, G])$, which in turn implies 
\begin{equation}\label{3eqn2}
|G/\Z(G)| \le |exp(\gamma_2(G))|^d.
\end{equation}
\vspace{.1in}

\noindent{\bf Problem 4.}  Classify all finite $p$-groups $G$ of nilpotency class $2$ for which equality holds in \eqref{3eqn2}.
\vspace{.1in}

Now we discuss some examples of infinite groups with finite central quotient. The most obvious example is the infinite cyclic group. Other obvious examples are the groups $G = H \times \mathcal{Z}$, where $H$ is  any finite group and  $\mathcal{Z}$ is the infinite cyclic group.  
Non-obvious examples include finitely generated  $FC$-groups, in which conjugacy class sizes are bounded, and certain Cernikov groups. We provide explicit examples in each case. Let $F_n$ be the free group on $n$ symbols and $p$ be a prime integer. Then the factor group $F_n/(\gamma_2(F_n)^p\gamma_3(F_n))$ is the required group of the first type, where $\gamma_2(F_n)^p = \gen{u^p \mid u \in \gamma_2(F_n)}$. Now let $H = \mathcal{Z}(p^{\infty}) \times A$ be the direct product of quasi-cyclic (Pr\"ufer) group $\mathcal{Z}(p^{\infty})$ and the cyclic group $A = \gen{a}$ of order $p$, where $p$ is a prime integer. Now consider the group $G = H \rtimes B$, the semidirect product of $H$ and  the cyclic group $B = \gen{b}$ of order $p$ with the action  by $x^b = x$ for all $x \in \mathcal{Z}(p^{\infty})$ and $a^b = ac$, where $c$ is the unique element of order $p$ in $\mathcal{Z}(p^{\infty})$. This group is suggested by V. Romankov and Rahul D. Kitture through ResearchGate, and is a Cernikov group. 

The following problem was suggested by R. Baer in \cite{rB38}.
\vspace{.1in}

\noindent{\bf Problem 5.} Let $A$ be an abelian group and $Q$ be a group. Obtain necessary and sufficient conditions on A and $G$ so that there exists a group $Q$ with $A \cong \Z(G)$ and $Q \cong G/A$.
\vspace{.1in}

This problem was solved by Baer himself for an arbitrary abelian group A and finitely generated abelian group G. Moskalenko \cite{aM68} solved this problem for an arbitrary abelian group A and  a periodic abelian group G. He \cite{aM71} also
solved this problem for arbitrary abelian group A and  a non-periodic abelian group G such that  the rank of $G/t(G)$ is $1$, where $t(G)$ denotes the tortion subgroup of $G$. If this rank is more than $1$, then he solved the problem when  $A$ is a torsion abelian group. For a given group $Q$, the existence of a group $G$ such that $Q \cong G/\Z(G)$ has been studied extensively under the theme `Capable groups'. However, to the best of our knowledge, Problem 5 has been poorly studied in full generality. Let us restate a special case of this problem in a little different setup. A pair  of groups $(A, Q)$, where  $A$ is an arbitray abelian group and $Q$ is an arbitrary group, is said to be a \emph{capable pair}  if there exists a  group $G$ such that $A \cong \Z(G)$ and $Q \cong G/\Z(G)$.  So, in our situation, the following problem is very interesting.
\vspace{.1in}

\noindent{\bf Problem 6.} Classify capable pairs $(A, Q)$, where  $A$ is an infinite abelian group and $Q$ is a finite group. 
\vspace{.1in}

Finally let us  get back to the situation when $G$ is a group with finite $\gamma_2(G)$ but infinite $G/\Z(G)$. The well known examples of such type are infinite extraspecial $p$-groups. Other class of examples can be obtained by taking a central product (amalgamated at their centers) of infinitely many copies of a $2$-generated finite $p$-group of class $2$ such that $\gamma_2(H) = \Z(H)$ is cyclic of order $q$, where $q$ is some power of $p$.  Notice that both of these classes consist of groups of nilpotency class $2$. Now if we take $G = X \times H$, where $X$ is an arbitrary group with finite $\gamma_2(X)$ and  $H$ is  with finite $\gamma_2(H)$ and infinite $H/\Z(H)$, then  $\gamma_2(G)$ is finite but  $G/\Z(G)$ is infinite. So we can construct nilpotent groups of arbitrary class and even non-nilpotent group with infinite central quotient and finite commutator subgroup. 

A non-nilpotent group $G$ is said to be \emph{purely non-nilpotent} if it does not have any non-trivial nilpotent subgroup as a direct factor. With the help of Rahul D. Kitture, we have also been able to construct purely non-nilpotent groups $G$ such that $\gamma_2(G)$ is finite but $G/\Z(G)$ is infinite. Let $H$ be an infinite extra-special $p$-group. Then we can always find a field $\mathcal{F}_q$, where $q$ is some power of a prime, containing all $p$th roots of unity. Now let $K$ be the special linear group $sl(p,  \mathcal{F}_q)$, which is a non-nilpotent group having  a central subgroup of order $p$. Now consider the group $G$ which is a central product of $H$ and $K$ amalgamated at $\Z(H)$. Then $G$ is a purely non-nilpotent group with the required conditions. It will be interesting to see more examples of this type which do not occur as a central product of such infinite groups of  nilpotency class $2$ with non-nilpotent groups.

By Proposition \ref{prop} we know that for an arbitrary group $G$ with finite $\gamma_2(G)$ but infinite  $G/\Z(G)$, the group   $G/\Z(G)$ has an infinite abelian group as a direct factor.  Further structural information is highly welcome.
\vspace{.1in}

\noindent{\bf Problem 7.}  Provide structural information of  the  group $G$ such that $\gamma_2(G)$ is finite but $G/\Z(G)$ is infinite?

\end{document}